\RequirePackage{fix-cm}
\documentclass[smallextended]{svjour3}       
\smartqed  
\usepackage{graphicx}

\usepackage{amssymb,amsmath,euscript}

\def\R{\mathbb{R}}
\def\F{\mathcal{F}}
\def\T{\mathcal{T}}

\def\A{\mathcal{A}}

\def\B{\mathcal{B}}

\def\I{\mathcal{I}}

\begin{document}

\title{The Sparsest Solutions to $Z$-Tensor Complementarity Problems \thanks{This research was supported by the National Natural Science Foundation of China (11301022,11431002), the State Key Laboratory of Rail Traffic Control and Safety, Beijing Jiaotong University (RCS2014ZT20, RCS2014ZZ01), and the Hong Kong Research Grant Council (Grant No.
PolyU 502111, 501212, 501913 and 15302114).}
}


\author{Ziyan Luo         \and
       Liqun Qi \and
        Naihua Xiu
}


\institute{Ziyan Luo \at
              The State Key Laboratory of Rail Traffic Control and Safety, Beijing
Jiaotong University, Beijing 100044, P.R. China; \\
              \email{starkeynature@hotmail.com}           
           \and
          Liqun Qi  \at
               Department of Applied Mathematics, The Hong Kong Polytechnic University, Hung Hom, Kowloon, Hong Kong, P.R.China;\\
               \email{liqun.qi@polyu.edu.hk}
          \and
          Naihua Xiu  \at
          Department of Mathematical, School of Sciences, Beijing Jiaotong University, Beijing, P.R.China; \\
          \email{nhxiu@bjtu.edu.cn}
}

\date{Received: date / Accepted: date}

\maketitle

\begin{abstract}
Finding the sparsest solutions to a tensor complementarity problem is generally NP-hard due to the nonconvexity and noncontinuity of the involved $\ell_0$ norm. In this paper, a special type of tensor complementarity problems with $Z$-tensors has been considered. Under some mild conditions, we show that to pursuit the sparsest solutions is equivalent to solving polynomial programming with a linear objective function. The involved conditions guarantee the desired exact relaxation and also allow to achieve a global optimal solution to the relaxed nonconvex polynomial programming problem. Particularly, in comparison to existing exact relaxation conditions, such as RIP-type ones, our proposed conditions are easy to verify.
\end{abstract}

\section{Introduction}
\label{intro}
The classical compressed sensing theory (see the pioneering work in \cite{CT2005,CRT2006,D2006}) has been gradually generalized and the nonlinear compressed sensing theory has attracted more and more attentions inspired by significant real-life applications such as sparse phase retrieval and sparse pulsation frequency detection in Asteroseimology (see \cite{B2013} and references therein). Among those nonlinear measurements, the polynomial structure has been employed in many applications cases, such as quadratic measurements in sparse signal recovery \cite{LV2013}, and nonlinear compressed sensing with polynomial measurements \cite{LO2014}. Besides the nonlinearity structure carried in many real-life applications, some priori information such as the nonnegativity, can be extensively encountered in communications, DNA microarrays, spectroscopy, tomography, network monitoring, and hidden Markov models \cite{DT2010,HMW2010,KDXH2011,LDL2010,SSMB2007}. In this regard, the optimality condition such as the KKT condition, which has been widely used in the optimization community, might be a good surrogate for the involved feasible set. Under some constraint qualifications, the original problem turns out to finding the sparsest solutions to a special nonlinear complementarity problem with polynomial structures. Mathematically, it can be formulated as
$$\quad{(PCP)~~~~}\\
\begin{array}{ll}
\min & \|x\|_0\\
{\rm s.t.} & F(x)\geq 0, x\geq 0,  \langle x, F(x)\rangle =0,
\end{array}$$
where $F=(f_1,\ldots, f_n)$ and all $f_i$'s are polynomial functions.

Apparently, when each $f_i$ in the aforementioned model $(PCP)$ is reduced to be affine (i.e., a linear function added by some constant), $(PCP)$ is exactly to find the sparsest solutions to a linear complementarity problem (LCP for short). The sparsest solution to LCP has been studied in \cite{CX2014,SZX2014}. However, to our best knowledge, related work on this topic are very limited, partially but essentially due to the complexity caused by the involved non-convex discontinuous objective function, and partially from the nonlinearity generated by the underlying complementarity constraint. For the former difficulty, many relaxation strategies have been explored by using different surrogates for the $\ell_0$ norm, such as the convex $\ell_1$ norm \cite{CRT2006,CT2005,D2006}, the non-convex $\ell_p$ norm \cite{CNZ2010,FL2009,XCXZ2012}, the reweighted $\ell_1$ norm \cite{CWB2008,CZ2014}, and so on. A natural but essential question arises: is it possible to get an exact solution of the original $\ell_0$ norm minimization problem by the relaxation counterpart? If so, what kind of properties should the involved data possess? For the linear measurement case, the well-known restricted isometry property (RIP for short) was introduced to guarantee the desired exactness, which has given a great explanation of the popularity of all sorts of random compressed sensing approaches \cite{C2008,CRT2006}. There are some other exact relaxation conditions on the coefficient matrix for linear constraints, such as the null space property \cite{Z2013}, the range space property \cite{Zhao2013}, the $s$-goodness property \cite{JN2011} and so on. Most of these exact relaxation properties are somehow not easy to verify. In \cite{Luo2014}, a generalized $Z$-matrix was introduced. Together with the nonnegativity of the right-hand side observation vector, it serves as an easy-to-check condition for the desired exact relaxation for the linear sparse optimization. For the nonlinear sparse optimization, such as the sparse LCP, the $Z$-matrix has been employed to guarantee the exact relaxation \cite{SZX2014}. Then how about more general polynomial cases? Can we find some $Z$-type condition to ensure an exact solution from the relaxation problem for the original sparse polynomial complementarity problem? This is our primary goal in this paper.

Recently, tensors, as a higher-order generalization of matrices, have been extensively studied \cite{CGL2008,LDV2000,Q2005,YY2010}, which is closely related to polynomials in terms of those coefficients. This allows us to write a polynomial equation system in a condense form with tensors. In this regard, when the constraint function in the aforementioned problem $(PCP)$ takes the form $F(x)=H(x)-b$ with some homogeneous polynomial function $H(x)$ and some vector $b\in \R^n$, then the feasible set can be reformulated as
$$\quad{(TCP)~~~~}\\ x\geq 0,~\A x^{m-1}-b \geq 0,~\langle x, \A x^{m-1}-b \rangle=0,$$
where $m-1$ is the degree of $H(x)$, $x^{m-1}$ is a rank one tensor of order $m-1$ and dimension $n$ with its $(i_1,\ldots, i_{m-1})$th entry $x_{i_1}\cdots x_{i_{m-1}}$, and $\A$ is an $m$th order $n$-dimensional tensor consisting of all the coefficients of $H(x)$ by means of $H(x)=\A x^{m-1}$. Here the tensor product $\A x^{m-1}$ is defined as $(\A x^{m-1})_i=\sum\limits_{i_2,\ldots,i_m=1}^na_{ii_2\ldots i_m}x_{i_2}\cdots x_{i_{m}}$, for all $i=1,\ldots, n$. Similarly, we can define $\A x^{m-k}$ as
$$(\A x^{m-k})_{i_1\ldots i_k}=\sum\limits_{i_{k+1},\ldots,i_m=1}^n a_{i_1\ldots i_ki_{k+1}\ldots i_m}x_{i_{k+1}}\cdots x_{i_{m}},~\forall i_1,\ldots, i_k=1,\ldots, n.\eqno(1.1)$$
The above $(TCP)$ is the so-called tensor complementarity problem which has been studied in \cite{CQW2015,SQ2014,SQi2014}. In this paper, we will focus on finding the sparsest solutions to a tensor complementarity problem which can be modeled as
$$\quad{(P_0)~~~~}\\
\begin{array}{ll}
\min & \|x\|_0\\
{\rm s.t.} & \A x^{m-1}-b\geq 0, x\geq 0,  \langle x, \A x^{m-1}-b\rangle =0.
\end{array}$$
Mathematically, problem $(P_0)$ is generally NP-hard due to the objective function $\|x\|_0$. Inspired by the scheme of the most popular convex relaxation, we could get the a linear surrogate $e^{T}x$ resulting from the nonnegativity constraint of $x$. But the nonlinearity from the tensor complementarity constraints cannot be easily handled and the existing exact relaxation conditions are not that appropriate since most of them are customized for linear systems.

In this paper, by employing $Z$-tensors (see e.g., \cite{DQW2013})and the least element theory in nonlinear complementarity problems \cite{T1974}, we present that if the involved $b$ is nonnegative and $\A$ is a $Z$-tensor, then a sparsest solution of the tensor complementarity problem can be achieved by solving the following polynomial programming problem:
$$\quad{(P_1)~~~~}\\
\begin{array}{ll}
\min & e^{T}x\\
{\rm s.t.} & \A x^{m-1}-b= 0, x\geq 0,~~~~~~~~~~~~~~~~~~~~~~~~~
\end{array}$$
where $e$ is the all one vector. In comparison to those existing exact relaxation conditions for general nonlinear sparse optimization problems \cite{B2013}, our conditions on the coefficients of the polynomial functions are easy to check. This is the main contribution of this paper. 

The rest of the paper is organized as follows.  The concepts of $Z$-tensor and $Z$-function are recalled and some useful properties are presented in Section 2. The $Z$-tensor complementarity problem is introduced and discussed in Section 3. The exact relaxation theorem is established in Section 4. Concluding remarks are drawn in Section 5.

For convenience of presentation, the following notations will be used throughout the paper. We use $\R^n$ and $\R^n_+$ to denote the $n$-dimensional Euclidean space and its nonnegative orthant respectively. $\R^{n\times n}$ is used to denote the space of all real $n\times n$ matrices. $\T_{m,n}$ is used to stand for the set of all real tensors with order $m$ and dimension $n$. Vectors are denoted by lowercase letters such as $x$, matrices are written as capital letters such as $A$, and tensors are written as calligraphic capital letters such as $\A$.

\section{$Z$-Tensors and $Z$-Functions}
\label{sec:1}
As a nonlinear generalization of $Z$-matrices with non-positive off-diagonal elements, the concept of off-diagonally antitone functions was first introduced by Rheinboldt in \cite{R1970}, which further leads to the definition of $Z$-functions as stated in \cite{T1974}. In \cite{I1992}, Isac has redefined the $Z$-functions equivalently by means of an implication system, which has been also widely used in the community of complementarity problems. In this section, the definition of $Z$-functions and some useful properties will be recalled. Particularly, to explore the $Z$-property of homogeneous polynomial functions, the $Z$-tensor and the partially $Z$-tensor will be introduced and analyzed.

\begin{definition}[Definition 3.2, \cite{I1992}]\label{def2} A mapping $F:\R^n\rightarrow \R^n$ is said to be a $Z$-function if for every $x$, $y$, $z\in\R^n_+$ such that $\langle x, y-z\rangle =0$, we have $\langle x, F(y)-F(z)\rangle\leq 0$.
\end{definition}



\begin{proposition}[Proposition 3.2, \cite{I1992}]\label{gradient} A G$\hat{a}$teaux continuous differentiable function $F:\R^n\rightarrow \R^n$ is a $Z$-function if and only if for any $x\in\R^n_+$, $\nabla F(x)$ is a $Z$-matrix.
\end{proposition}

\begin{lemma}[\cite{I1992}]\label{implication} If $F:\R^n\rightarrow \R^n$ is a $Z$-function, then the following implication holds:
$$x\in\R^n_+, y\in\R^n_+, \langle x,y\rangle=0\Rightarrow \langle x, F(y)-F(0)\rangle \leq 0.\eqno(2.1)$$
Moreover, if $F(x)=Ax$ is a linear function, then $A$ is a $Z$-matrix, which is equivalent to the following implication:
$$x\in\R^n_+, y\in\R^n_+, \langle x,y\rangle=0\Rightarrow \langle x, Ay\rangle \leq 0.\eqno(2.2)$$
\end{lemma}

\vskip 2mm

It is known that for any matrix $A$, $Ax$ is a $Z$-function if and only if $A$ is a $Z$ matrix, i.e., all of its off-diagonal entries are non-positive. This concept has been extended to the higher order tensors as stated below.


\begin{definition}[\cite{ZQZ2014}]\label{def1} Let $\A=(a_{i_1\ldots i_m})\in \T_{m,n}$. $\A$ is called a $Z$-tensor if all its off-diagonal entries are nonpositive, i.e., $a_{i_1\ldots i_m}\leq 0$ when $\delta_{i_1,\ldots, i_m}=0$.
\end{definition}

Another concept called the partially $Z$-tensor is introduced here.\vskip 2mm

\begin{definition}\label{def3} Let $\A=(a_{i_1\ldots i_m})\in \T_{m,n}$. We call $\A$ a partially $Z$-tensor if for any $i_1\in [n]$, $a_{i_1i_2\ldots i_m}\leq 0$ for all $i_2,\ldots, i_m$ satisfying $i_1\notin \{i_2,\ldots, i_m\}$.
\end{definition}

Obviously, a $Z$-tensor is a partially $Z$-tensor, and both $Z$- and partially $Z$-tensors of order $m=2$ are exactly $Z$-matrices. Thus, these two concepts both can be regarded as extensions of the $Z$-matrix. Properties on these two types of tensors are discussed as follows which will play an essential role in the sequel analysis. \vskip 2mm

\begin{theorem}\label{partial} For any given $\A\in \T_{m,n}$, we have

(i) if $\A$ is a partially $Z$-tensor, then the implication (2.1) holds with $F(x)=\A x^{m-1}$;

(ii) if $\A$ is a $Z$-tensor, then $F(x)=\A x^{m-1}$ is a $Z$-function.
\end{theorem}

\noindent{\bf Proof.} (i) Suppose $x$, $y\in\R^n_+$ with $\langle x,y\rangle =0$. Easily we can get
$$x_i\geq 0, y_i\geq 0, x_iy_i=0,\forall i\in[n].\eqno(2.3)$$
Thus,  \begin{eqnarray}
       & &\langle x, F(y)-F(0)\rangle\nonumber\\
       &=& \langle x, \A y^{m-1}\rangle=\sum\limits_{i=1}^n x_i\left(\A y^{m-1}\right)_i\nonumber\\
       &=& \sum_{i=1}^n x_i \sum_{i_2,\ldots, i_m=1}^n a_{ii_2\ldots i_m}y_{i_2}\cdots y_{i_m} \nonumber \\
         &=&  \sum_{i=1}^n \left(\sum_{\begin{subarray}{c}
         i_2,\ldots, i_m=1 \\
         i \notin \{i_2,\ldots, i_m\}
         \end{subarray}}^n a_{ii_2\ldots i_m} x_i y_{i_2}\cdots y_{i_m}+\sum_{\begin{subarray}{c}
         i_2,\ldots, i_m=1 \\
         i \in \{i_2,\ldots, i_m\}
         \end{subarray}}^n  a_{ii_2\ldots i_m} x_i y_{i_2}\cdots y_{i_m}\right) \nonumber \\
          &=&  \sum_{i=1}^n \sum_{\begin{subarray}{c}
         i_2,\ldots, i_m=1 \\
         i \notin \{i_2,\ldots, i_m\}
         \end{subarray}}^n  a_{ii_2\ldots i_m} x_i y_{i_2}\cdots y_{i_m} \nonumber\\
          &\leq & 0, \nonumber
      \end{eqnarray}
where the last equality follows from (2.3) and the last inequality follows from Definition \ref{def3}.

\noindent (ii) Invoking Proposition \ref{gradient}, it suffices to show that for any $x\in\R^n_+$, $\nabla_x\left(\A x^{m-1}\right)$ is a $Z$-matrix. Combining with Lemma \ref{implication}, we only need to show the implication (2.2) holds with $A:=\A x^{m-1}$ for any given $x\in\R^n_+$. Let $y$, $z$ be any two nonnegative vectors, and $\langle y, z\rangle=0$. It yields that
\begin{eqnarray}
  &&\langle y, \nabla_x\left(\A x^{m-2}\right)z\rangle \nonumber\\
  &=& \sum_{i=1}^ny_i\sum_{i_2,\ldots, i_m=1}^n \left(a_{ii_2\ldots i_m}+a_{i_2 i\ldots i_m}+\ldots+a_{i_2\ldots i_mi}\right)x_{i_2}\cdots x_{i_{m-1}}z_{i_m} \nonumber \\
   &=&  \sum_{i=1}^n\left(\sum_{\begin{subarray}{c}
   i_2,\ldots, i_m=1 \\
   i_m\neq i
   \end{subarray}}^n \left(a_{ii_2\ldots i_m}+a_{i_2 i\ldots i_m}+\ldots+a_{i_2\ldots i_mi}\right)x_{i_2}\cdots x_{i_{m-1}}z_{i_m}y_i\right) \nonumber \\
   &\leq &0, \nonumber
\end{eqnarray}
where the second equality is from the complementarity of $y$ and $z$, and the last inequality is from the fact that $\A$ is a $Z$-tensor. Thus $\F_{\A}$ is a $Z$-function. \qed

\section{$Z$-Tensor Complementarity Problems}
\label{sec:2}
It is known that a tensor complementarity problem always takes the form
$$\quad{\left(TCP(\A,b)\right)~~~~}\\
x\geq 0, ~\A x^{m-1}-b\geq 0,~ \langle x,~\A x^{m-1}-b\rangle=0,
$$
\noindent which is actually a special nonlinear complementarity problem. When the involved tensor $\A$ is a $Z$-tensor, the corresponding $(TCP(\A,b))$ is called a $Z$-tensor complementarity problem. In this section, we will concentrate on exploiting the properties of such a special class of tensor complementarity problems. We start with recalling a nice property possessed by general nonlinear complementarity problems with $Z$-functions.

\begin{theorem}[Ex. 3.7.21, \cite{FP2003}]\label{least} Let $F:\R^n\rightarrow \R^n$ be a continuous $Z$-function. Suppose that the following nonlinear complementarity problem
$$ \quad{(NCP(F))~~~~}\\
 x\geq 0, F(x)\geq 0, \langle x,F(x)\rangle=0$$
is feasible, i.e., $\mathcal{F}:=\{x\in\R^n: x\geq 0, F(x)\geq 0\}\neq \emptyset$. Then ${\mathcal{F}}$ has a unique least element $x^*$ which is also a solution to $(NCP(F))$.
\end{theorem}

Inspired by the relation between $Z$-tensors and $Z$-functions, we can easily get the following properties for $Z$-tensor complementarity problems.

\begin{corollary}\label{eq1} Let $\A$ be a $Z$-tensor and $b\in\R^n$. Suppose that the tensor complementarity problem $(TCP(\A,b))$ is feasible, i.e., $\mathcal{F}:=\{x\in\R^n: x\geq 0, \A x^{m-1}-b\geq 0\}\neq \emptyset$. Then ${\mathcal{F}}$ has a unique least element $x^*$ which is also a solution to $(TCP(\A, b))$.
\end{corollary}
\noindent{\bf Proof.} Theorem \ref{partial} tells us that $\A x^{m-1}$ is a $Z$-function. Utilizing Proposition \ref{gradient}, it is easy to verify that $\A x^{m-1}-b$ is also a $Z$-function for any $b\in\R^n$. Thus, the desired result follows directly from Theorem \ref{least}. \qed

With a nonnegative $b$ and a partially $Z$-tensor $\A$ in $(TCP(\A,b))$, the tensor complementarity problem can be equivalent to a multi-linear equation with nonnegative constraints.

\begin{proposition}\label{equi1} Let $\A$ be a partially $Z$-tensor and $b\in\R^n_+$. The following two systems are equivalent:

(i) $x\in\R^n_+$, $\A x^{m-1}-b\in\R^n_+$, $\langle x, \A x^{m-1}-b\rangle=0$;

(ii) $x\in\R^n_+$, $\A x^{m-1}-b=0$.
\end{proposition}
\noindent{\bf Proof.} Trivially, any solution to system (ii) is a solution to system (i). Let $y$ be any solution to system (i). Since $\A$ is a partially $Z$-tensor and $b\in\R^n_+$, invoking Theorem \ref{least}, it yields that
\begin{eqnarray}
0 &\geq & \langle \A y^{m-1}-b, \A y^{m-1}\rangle \nonumber\\
&=& \langle \A y^{m-1}-b,\A y^{m-1}-b\rangle +\langle \A y^{m-1}-b,b\rangle \nonumber\\
&\geq & \|\A y^{m-1}-b\|_2^2. \nonumber
\end{eqnarray}
This indicates that $\A y^{m-1}-b=0$, which implies that $y$ is a solution to (ii). \qed

Note that $Z$-tensors are partially $Z$-tensors. Thus the results in the above proposition hold for $Z$-tensors.\vskip 2mm

\begin{corollary}\label{co1} Let $\A$ be a $Z$-tensor and $b\in\R^n_+$. The following two systems are equivalent:

(i) $x\in\R^n_+$, $\A x^{m-1}-b\in\R^n_+$, $\langle x, \A x^{m-1}-b\rangle= 0$;

(ii) $x\in\R^n_+$, $\A x^{m-1}-b=0$.
\end{corollary}

Utilizing the aforementioned equivalence, we can characterize the feasibility of $TCP(\A,b)$ in terms of the consistency of the corresponding nonnegative constrained multi-linear equation. Before stating the feasibility, we recall the definition of $M$-tensors, which form an important subclass of $Z$-tensors. \vskip 2mm

\begin{definition}[\cite{ZQZ2014}] \label{def5} Let $\A\in\T_{m,n}$ be a $Z$-tensor with $\A=s \I-{\mathcal{B}}$, where $\I$ is the identity tensor whose diagonal entries are $1$ and others $0$, $\mathcal{B}$ is a nonnegative tensor and $s\in\R^n_+$. If $s\geq \rho({\mathcal{B}})$, then $\A$ is called an $M$-tensor. If $s>\rho({\mathcal{B}})$, then $\A$ is called a strong $M$-tensor. Here $\rho({\mathcal{B}})$ stands for the spectral radius of $\mathcal{B}$.
\end{definition}

\begin{proposition}\label{feasibility} If $\A$ is a strong $M$-tensor and $b\in\R^n_+$, then $(TCP(\A,b))$ is feasible.
\end{proposition}
\noindent{\bf Proof.} Let $\A=s\I-\B$ be a strong $M$-tensor with $\B\geq 0$ and $s>\rho (\B)$. Invoking the equivalence as established in Corollary \ref{co1}, it suffices to show that that there exists some nonnegative $x$ such that $\A x^{m-1}=b$. Let $T_{s,\B,b}: \R^n_+\rightarrow \R^n_+$ be the mapping defined as follows:
$$T_{s,\B,b}(x):=\left(s^{-1}\B x^{m-1}+s^{-1}b\right)^{[\frac{1}{m-1}]},~~\forall x\in\R^n,$$
where $x^{[\frac{1}{m-1}]}$ is the vector with its $i$th component $x_i^{[\frac{1}{m-1}]}$, for all $i=1,~\ldots,~n$. Easily, we can find that the required nonnegative solution $x$ is exactly a fixed point of this mapping $T_{s,\B,b}$. Besides, since $\A$ is a strong $M$-tensor, applying Theorem 3 in \cite{DQW2013}, there always exists a positive $z$ such that $\A z^{m-1}>0$. Denote
$$\alpha:=\min\limits_{i\in [n]}\frac{b_i}{\left(\A z^{m-1}\right)_i},~\textrm{and}~\beta:=\max\limits_{i\in [n]}\frac{b_i}{\left(\A z^{m-1}\right)_i}.$$
Obviously, $$0\leq \alpha \A z^{m-1}\leq b\leq \beta \A z^{m-1}.$$
Set $v:=\A \left(\alpha^{\frac{1}{m-1}}z\right)^{m-1}$ and $w:=\A \left(\beta^{\frac{1}{m-1}}z\right)^{m-1}$. Therefore,
$$\alpha^{\frac{1}{m-1}}z=T_{s,\B,v}\left(\alpha^{\frac{1}{m-1}}z\right)\leq T_{s,\B,b}\left(\alpha^{\frac{1}{m-1}}z\right),$$ $$\beta^{\frac{1}{m-1}}z=T_{s,\B,w}\left(\alpha^{\frac{1}{m-1}}z\right)\geq T_{s,\B,b}\left(\beta^{\frac{1}{m-1}}z\right)$$
Note that $T_{s,\B,b}$ is an increasing continuous mapping on $\R^n_+$. By employing the fixed point theorem in \cite{A1976} (also see Theorem 3.1 in \cite{DW2015}), there exists at least one fixed point $x$ of $T_{s,\B,b}$ such that $0\leq \alpha^{\frac{1}{m-1}}z\leq x\leq \beta^{\frac{1}{m-1}}z$. This completes the proof.\qed

\section{Exact Relaxation}
\label{sec:3}

Now we are in a position to establish the exact relaxation theorem for the $\ell_0$ norm minimization problem $(P_0)$.\vskip 2mm

\begin{theorem}\label{exact-relaxation} Let $\A$ be a $Z$-tensor and $b\in\R^n_+$. If the problem $(P_0)$ is feasible, then $(P_0)$ has a solution $x^*$ which is also the unique solution to the problem $(P_1)$.
\end{theorem}
\noindent{\bf Proof.} Invoking Theorem \ref{least}, we know that $\A x^{m-1}$ is a $Z$-function. Together with Proposition \ref{gradient}, it is easy to verify that $\A x^{m-1}-b$ is also a $Z$-function. Thus, Theorem \ref{least} tells us that there exists a unique least element $x^*$ in ${\mathcal{F}}$ as defined in Corollary \ref{eq1}, which is also a solution of the tensor complementarity problem. The nonnegativity constraint directly yields that $x^*$ is one of the sparsest solutions of $(P_0)$. Utilizing the equivalence as shown in Corollary \ref{co1}, $x^*$ is definitely the unique solution of $(P_1)$ by the fact that it should be the least element in $\{ x\in\R^n: x\geq 0,\A x^{m-1}=b\}$. This completes the proof.\qed

\begin{corollary} Let $\A$ be a strong $M$-tensor and $b\in\R^n_+$. Then problem $(P_1)$ is uniquely solvable and the unique solution is also an optimal solution to problem $(P_0)$.
\end{corollary}
\noindent{\bf Proof.} This follows directly from Proposition \ref{feasibility} and Theorem \ref{exact-relaxation}.\qed

Some extended result on exact relaxation theory is discussed as follows. For any matrix $P=\left(p_{ij}\right)\in\R^{n\times n}$, we define a linear operator $P\A: \T_{m,n}\rightarrow \T_{m,n}$ as follows:
 $$ \left(P \A  \right)_{i_1i_2\ldots i_m} =\sum\limits_{i=1}^n p_{i_1 i}a_{ii_2\ldots i_m},~\forall \A=\left(a_{i_1i_2\ldots i_m}\right)\in\T_{m,n}.$$
This is also treated as a special tensor-matrix product in the tensor community.(see e.g., \cite{KB2009}) Evidently, if $P$ is an invertible matrix, then $P\T_{m,n}=\T_{m,n}$. However, this operator cannot preserve the $Z$-property for $Z$-tensors. Note that the right-hand side of the multi-linear system $\A x^{m-1}=b$ is actually a condense form of $n$ homogeneous polynomials of degree $m$ with any row tensor $\A_i$ formed by the corresponding coefficients of the $i^{\textrm{th}}$ polynomial. If we change the order of these equations in the multi-linear system, the solution set will not be affected. This observation allows us to generalize the result in Theorem \ref{exact-relaxation}. For convenience, we use ${\mathbb{P}}^Z_{m,n}$ to denote the set of all tensors of order $m$ and dimension $n$ which can be transformed to $Z$-tensors with some permutation matrices, i.e.,
$${\mathbb{P}}^Z_{m,n}:=\{\A\in \T_{m,n}: P\A \textrm{~is~a~} Z\textrm{-tensor},P \textrm{~is~a~permutation~matrix} \}$$

\begin{corollary}\label{exact2} Suppose $\A \in {\mathbb{P}}^Z_{m,n}$ and $b\in\R^n_+$. If the problem $(P_0)$ is feasible, then $(P_0)$ has a unique solution $x^*$ which is also the unique solution to the problem $(P)$.
\end{corollary}


\section{Conclusions}
To pursuit the sparsest solutions to a tensor complementarity problem can be formulated as an $\ell_0$ norm minimization with tensor complementarity constraints, which is always NP-hard. Based on the properties of $Z$-tensors, we show that one of the sparsest solutions of the $Z$-tensor complementarity problem can be achieved in polynomial time by solving a polynomial programming problem with a linear objective function, such as the Gauss-Seidel iteration method proposed in \cite{DW2015}. The involved condition on the input of data is indeed an exact relaxation condition for the original $\ell_0$ norm minimization. It is worth mentioning that in comparison to other existing exact relaxation conditions in the community of sparse optimization or compressed sensing, our proposed condition is easy to verify. Including the sparse tensor complementarity problem as a special case, the topic of sparse optimization with general nonlinear complementarity constraints deserves further study.

\bibliographystyle{spmpsci}

\bibliography{ZTCP}



%
%

\end{document}